# Primes In Cubic Arithmetic Progressions

N. A. Carella, February, 2013.

*Abstract*: This work proposes a proof of the simplest cubic primes counting problem. It shows that the subset of primes $P_f = \{ p = n^3 + 2 \text{ is prime} : n \in \mathbb{Z} \}$ is an infinite subset of primes. Further, the expected order of magnitude of the cubic primes counting function is established.

*Mathematics Subject Classifications*: 11A41, 11N05, 11N32, 11N13.

Keywords: Distribution of Primes, Cubic Primes, Bouniakowsky Conjecture, Prime Diophantine Equations.

## 1. Introduction

The *cubic primes conjecture* claims that certain Diophantine equations $y = ax^3 + bx^2 + cx + d$ have infinitely many prime solutions $y = p$ as the integer $x = n \in \mathbb{Z}$ varies. More generally, the *Bouniakowsky conjecture* claims that for any irreducible polynomial $f(x) \in \mathbb{Z}[x]$ over the integers of *fixed divisor* $\text{div}(f) = 1$, and degree $\deg(f) \geq 1$, the Diophantine equation $y = f(x)$ has infinitely many prime solutions $y = p$ as the integer $x = n \in \mathbb{Z}$ varies. The case of linear polynomials $f(x) = ax + b$, $\gcd(a, b) = 1$, is known as Dirichlet Theorem, and the related case of linear polynomials $f(x) = ax + b$ with $k \geq 1$ consecutive primes for some $\gcd(a, b) = 1$, is known as the Green-Tao Theorem.

The fixed divisor $\text{div}(f) = \gcd(f(\mathbb{Z}))$ of a polynomial $f(x) \in \mathbb{Z}[x]$ over the integers is the greatest common divisor of its image $f(\mathbb{Z}) = \{ f(n) : n \in \mathbb{Z} \}$ over the integers. The fixed divisor $\text{div}(f) = 1$ if the congruence equation $f(x) \equiv 0 \bmod p$ has $w(p) < p$ solutions for all primes $p \leq \deg(f)$, see [FI, p. 395].

Some irreducible polynomials of fixed divisors $\text{div}(f) \neq 1$ can generate at most one prime. For example, the irreducible polynomials:

(i) $f_1(x) = x(x + 1) + 2$ has fixed divisor $\text{div}(f_1) = 2$,
(ii) $f_2(x) = x(x + 1)(x + 2) + 3$ has fixed divisor $\text{div}(f_2) = 3$, and
(iii) $f_3(x) = x^p - x + p$ has fixed divisor $\text{div}(f_3) = p$.

Since the fixed divisors div(f) > 1, these polynomials can generate up to one prime. But, the irreducible polynomials

(iv) $f_4(x) = x^2 + 1$ has fixed divisor div($f_4$) = 1,

(v) $f_5(x) = x^3 + 2$ has fixed divisor div($f_5$) = 1,

so these polynomials can generate infinitely many primes.

Reducible polynomials $f(x)$ of degrees $d = \deg(f)$ can generate up to $d \geq 1$ primes. For example, the reducible polynomial $f(x) = (g(x) \pm 1)(x^3 + 2)$ can generates up to $\deg(g) \geq 1$ primes for any polynomial $0 \neq g(x) \in \mathbb{Z}[x]$ of degree $\deg(g) \geq 1$.

**Cubic Primes Conjecture.** For a fixed noncubic integer $k \in \mathbb{Z}$, and a cubic polynomial $f(x) = x^3 + k \in \mathbb{Z}[x]$ of fixed divisor div($f(x)$) = 1, there are infinitely many primes $p = n^3 + k$ as $n$ varies over the integers. The number of cubic primes $p$ up to $x$ is asymptotically given by

$$\pi_f(x) = \frac{x^{1/3}}{\log x} \prod_{\substack{p \equiv 1 \bmod 3, \\ \gcd(p,k)=1}} \left(1 - \frac{2\chi(-k)}{p-1}\right) + O\left(\frac{x^{1/3}}{\log^2 x}\right), \qquad (1)$$

where the cubic symbol $\chi$ is defined by

$$\chi(-k) = \begin{cases} 1 & \text{if } -k \text{ is a cubic residue mod } p, \\ -1/2 & \text{if } -k \text{ is not a cubic residue mod } p. \end{cases} \qquad (2)$$

The noncubic integer $k$ condition preempts the reducible polynomials $x^3 + k = (x + a)(x^2 - ax + a)$ with $k = a^3$. This conjecture is discussed in [RN, p. 407], [NW, p. 348]; tables and computational works are given in [CW] and [SK] respectively. A result on the average order of the number of primes $p = x^3 + k$, for $k$ in certain range, but not for any single $k$, is given in [FZ]. The general conjecture and detailed discussions of the prime values of polynomials appear in [RN, p. 387], [GL, p. 17], [FI, p. 395], [NW, p. 405], [WK], [PP], and related topics in [BR], [GM], [MA], [PZ, p. 33], et alii.

## 2. The Proof of the Main Result

This Section provides the details of a possible proof of the simplest cubic primes counting problem. In particular, the cardinality of the subset of cubic primes $\{ p = n^3 + 2 : n \geq 1 \}$ is infinite.

The analysis of primes in the cubic arithmetic progression $\{ n^3 + 2 : n \geq 1 \}$, which is specified by one of the simplest cubic polynomial, draws on the standard results on the cubic residuacity of the integer 2, the properties of the subset of integers $D_f = \{ n \in \mathbb{N} : x^3 + 2 \equiv 0 \bmod n \text{ for some } x \in \mathbb{Z} \}$, and the corresponding Dirichlet series $D(s, f) = \sum_{n \in D_f} (\mu(n) \log n) n^{-s}$.

***Theorem* 1.** The Diophantine equation $y = x^3 + 2 \in \mathbb{Z}[x, y]$ has infinitely many primes $y = p$ solutions as $x = n \in \mathbb{Z}$ varies over the integers.

***Proof***: Since the congruence $x^3 + 2 \equiv 0 \bmod p$ has $< 3$ solutions for each prime $p \le 3 = \deg(f)$, the polynomial $f(x) = x^3 + 2$ has the fixed divisor $\operatorname{div}(f) = 1$, confer [FI, p. 395]. Now, select an appropriate weighted finite sum over the integers, see (8), and rewrite it as

$$\sum_{n^3 \le x} n\Lambda(n^3 + 2) = -\sum_{n^3 \le x} n \sum_{d \mid n^3+2} \mu(d) \log d$$

$$= -\sum_{d \le x+2} \mu(d) \log d \left( \sum_{n \le x^{1/3},\, n^3 \equiv -2 \bmod d} n \right). \tag{3}$$

This transformation follows from the identity $\Lambda(n) = -\sum_{d \mid n} \mu(d) \log d$, see Lemma 2, and inverting the order of summation. Other examples of inverting the order of summations are given in [MV, p. 35], [RE, p. 27], [RM, p. 216], [SP, p. 83], [TM, p. 36], et alii.

The congruence $n^3 \equiv -2 \bmod d$ is equivalent to the restrictions $p \mid d \Rightarrow p = 3m + 2$, or $p = u^2 + 27v^2$, or 3. These conditions specify the subset of integers $D_f = \{ d \in \mathbb{N} : x^3 + 2 \equiv 0 \bmod d \text{ for some } x \in \mathbb{Z} \}$, see Lemma 3. The solution counting function $\rho(d) \ge 0$ given below, is defined in (14).

Applying Lemma 4, (with $a = -2$ and $q = d$), using the fact that $\rho(d) \ge 1$ for $d \in D_f$, and simplifying yield

$$\sum_{n^3 \le x} n\Lambda(n^3 + 2) = -\sum_{d \le x+2} \mu(d) \log d \left( \frac{\rho(d)}{2d} \left(x^{1/3}\right)^2 + o\left(\frac{\rho(d)}{d} x^{2/3}\right) \right)$$

$$\ge -c_1 x^{2/3} \sum_{d \le x+2,\, d \in D_f} \frac{\mu(d) \log d}{d} + o\left( x^{2/3} \left| \sum_{d \le x+2,\, d \in D_f} \frac{\mu(d) \log d}{d} \right| \right), \tag{4}$$

where $c_1 > 0$ is a constant.

Applying Lemma 5 (for example, at $s = 1$, for $q = d$), the previous equation becomes

$$\sum_{n^3 \le x} n\Lambda(n^3 + 2) \ge c_1 x^{2/3} \left( c_2 + O\left(e^{-c(\log x)^{1/2}} \log x\right) \right) + o\left( x^{2/3} \left( c_2 + O(e^{-c(\log x)^{1/2}} \log x) \right) \right)$$

$$= c_3 x^{2/3} + O\left( x^{2/3} e^{-c(\log x)^{1/2}} \log x \right), \tag{5}$$

where $c_1, c_2, c_3 > 0$, and $c > 0$ are nonnegative constants. On the contrary, assume that there are finitely many cubic primes $p = n^3 + 2 < x_0$, where $x_0 > 0$ is a large constant. For example, $\Lambda(n^3 + 2) \ne \Lambda(p)$ with $p$ prime for all $n^3 + 2 > x_0$, $n \in \mathbb{N}$. Then

$$c_3 x^{2/3} + O\left(x^{2/3} e^{-c(\log x)^{1/2}} \log x\right) \leq \sum_{n^3 \leq x} n\Lambda(n^3 + 2)$$

$$= \sum_{n^3 \leq x_0} n\Lambda(n^3 + 2) + \sum_{x_0 < n^3 \leq x} n\Lambda(n^3 + 2) \qquad (6)$$

$$\leq x_0^{2/3} + c_4 x^{1/2} \log^2 x,$$

where $c_4 > 0$ is a constant, for all sufficiently large real numbers $x \geq 1$, see Lemma 6. But this is a contradiction for all sufficiently large real number $x > x_0$. ∎

Presumably, the analysis for pure cubic polynomials $f(x) = x^3 + k \in \mathbb{Z}[x]$ uses the same type algebraic and analytical manipulations as the simplest case $f(x) = x^3 + 2$, confer [BT, Theorem 7.1.7] for information on cubic residues. The general case $f(x) = ax^3 + bx^2 + cx + d \in \mathbb{Z}[x]$ probably can be handled via the weighted finite sum

$$\sum_{n^3 \leq x} n\Lambda(an^3 + bn^2 + cn + d) \quad \text{or} \quad \sum_{n^3 \leq x} n^2 \Lambda(an^3 + bn^2 + cn + d) \qquad (7)$$

or similar finite sums. Moreover, it will require concrete descriptions of the subset of integers $D_f = \{ d \in \mathbb{N} : f(x) \equiv 0 \bmod d \text{ for some } x \in \mathbb{Z} \}$ defined by admissible cubic polynomials, the solution counting function $\rho(d)$, the convergence properties of the corresponding Dirichlet series $D(s, f) = \sum_{n \in D_f} (\mu(n) \log n) n^{-s}$, and possibly other data on the associated cubic and quadratic fields.

Apparently, this analysis works well with any weighted finite sum of the form

$$\sum_{n \leq x} n^k \Lambda(n^3 + 2) = M_k(x) + E_k(x) \qquad (8)$$

with $k \geq 1$. But fails if $k = 0$ because the error term $E_0(x)$ is larger than the main term $M_0(x)$. It should also be noted that other weighted finite sums over the integers such as

$$\sum_{n \leq x} \varphi(n) \Lambda(n^3 + 2), \quad \sum_{n \leq x} \sigma(n) \Lambda(n^3 + 2), \quad \sum_{n \leq x} \tau(n) \Lambda(n^3 + 2), \qquad (9)$$

where $\varphi$, $\sigma$, and $\tau$ are the totient function, and the divisors functions respectively, yield precisely the same result. The analysis of these later finite sums are quite similar to the Titchmarsh divisor problem over cubic arithmetic progressions, see [CM, p. 172] for an introduction to this problem.

The topic of primes in quadratic, cubic, quartic arithmetic progressions, the Bouniakowsky conjecture, and in general the Bateman-Horn conjecture, is a rich area of research involving class fields theory and analytic number theory.

## 3. Elementary Foundation
The elementary underpinning of Theorem 1 is assembled here. The basic definitions of several number theoretical functions, and a handful of Lemmas are recorded here. The proofs of all these lemmas use elementary methods.

Let $n \in \mathbb{N} = \{0, 1, 2, 3, \ldots\}$ be a nonnegative integer. The Mobius function is defined by

$$\mu(n) = \begin{cases} (-1)^v & \text{if } n = p_1 p_2 \cdots p_v, \\ 0 & \text{if } n \neq \text{squarefree}. \end{cases} \qquad (10)$$

And the vonMangoldt function is defined by

$$\Lambda(n) = \begin{cases} \log p & \text{if } n = p^k, k \geq 1, \\ 0 & \text{if } n \neq p^k, k \geq 1. \end{cases} \qquad (11)$$

**Formula for the vonMangoldt Function**

***Lemma 2.*** Let $n \geq 1$ be an integer, and let $\Lambda(n)$ be the vonMangoldt function. Then

$$\Lambda(n) = -\sum_{d \mid n} \mu(d) \log d. \qquad (12)$$

***Proof***: Use Mobius inversion formula on the identity $\log n = \sum_{d \mid n} \Lambda(d)$ to verify this claim. ∎

Extensive details and other identities and approximations of the vonMangoldt are discussed in [FI], [HW], [HN], et alii.

**Powers Sums Over Cubic Arithmetic Progressions**

Let $a \geq 1$, and $q \geq 1$ be integers, $\gcd(a, q) = 1$. The element $a \in \mathbb{Z}_q = \{0, 1, 2, \ldots, q-1\}$ is called a *cubic residue* if the congruence $x^3 \equiv a \bmod q$ has a solution $x = b \in \mathbb{Z}_q$. Otherwise, $a$ is a *cubic nonresidue*. The cubic symbol is defined by

$$\left(\frac{a}{\pi}\right)_3 = \begin{cases} 1 & \text{if } a \text{ is a cubic residue}, \\ \omega^m & \text{if } a \text{ is a cubic nonresidue}, \\ 0 & \text{if } \gcd(a, \pi) \neq 1, \end{cases} \qquad (13)$$

where $\pi \in \mathbb{Z}[\omega]$ is a prime, and the complex number $\omega = (-1 + \sqrt{-3})/2$ is a cubic root of unity, and $m = 1, 2$, see [BT, p. 235].

Let $\rho(q) = \#\{x : f(x) \equiv 0 \bmod q\}$ denotes the solution counting function. The value $\rho(q) \geq 0$ depends on the primes divisors of $q$. It is a simple task to verify that every integer $a \geq 1$ is a cubic residue modulus a prime $p = 3m + 2$ or 3. Accordingly, the congruence $n^3 \equiv a \bmod p$ has a single solution $\rho(q) = 1$ if $p = 3m + 2$ or 3. On the other hand, for a prime $p = 3m + 1 > 3$, the congruence $n^3 \equiv a \bmod q$ has no solution $\rho(q) = 0$ or $\rho(q) = 3$ solutions.

**Lemma 3.** (Gauss) Let $p = 3n + 1$ be a prime number. Then
(i) The integer 2 is cubic residue modulus $p$ if and only if $p = u^2 + 27v^2$.
(ii) The integer 2 is cubic nonresidue modulus $p$ if and only if $p = 4u^2 + 2uv + 7v^2$.

The cubic residuacity of numbers and related concepts are explicated in [BT, p. 213], and [CX, p. 80]. These concepts come into play in the analysis of power sums over cubic arithmetic progressions.

The previous Lemma quickly leads to a formula for the solutions counting function

$$\rho(q) = \prod_{p | q} \rho(p), \quad \text{where} \quad \rho(p) = \begin{cases} 1 & \text{if } p = 3n + 2 \text{ or } 3, \\ 3 & \text{if } p = 3n + 1 = u^2 + 27v^2, \\ 0 & \text{if } p = 3n + 1 \neq u^2 + 27v^2, \end{cases} \quad (14)$$

for squarefree $q \geq 1$, which is a multiplicative function.

**Lemma 4.** Let $a \geq 1$, and $q \geq 1$ be integers, $\gcd(a, q) = 1$. Assume that $a$ is a cubic residue modulus $q$. Then

(i) $\displaystyle\sum_{n \leq x, \ n^3 \equiv a \bmod q} n = \frac{\rho(q)}{2q} x^2 + o(\frac{\rho(q)}{q} x^2),$ \qquad (ii) $\displaystyle\sum_{n \leq x, \ n^3 \equiv a \bmod q} n \geq \frac{1}{2q} x^2 + o(\frac{1}{q} x^2),$ \qquad (15)

where $\rho(q) \geq 0$ is a function of the prime divisors of $q$, for all sufficiently large real numbers $x \geq 1$.

**Proof** (ii): The assumption that $a \geq 1$ is a cubic residue immediately implies that the congruence $x^3 \equiv a \bmod q$ has at least one solution, and $\sum_{n \leq x, \ n^3 \equiv a \bmod q} n \geq 1$. Now, for each solution $x = b$ of the congruence $x^3 \equiv a \bmod q$ with $0 < b < q$, the sequence of integers in the arithmetic progression $n = qm + b$, with $0 \leq m \leq (x - b)/q$, are also solutions of the congruence. Thus, expanding, and simplifying the finite sum return

$$\sum_{n \leq x, \ n^3 \equiv a \bmod q} n = \sum_{n \leq x, \ n \equiv b \bmod q} n$$

$$= q \sum_{m \leq (x-b)/q} m + b \sum_{m \leq (x-b)/q} 1 \qquad (16)$$

$$= \frac{q}{2} \left[\frac{x-b}{q}\right]\left(\left[\frac{x-b}{q}\right] + 1\right) + b\left[\frac{x-b}{q}\right],$$

where $[z] = z - \{z\}$ is the largest integer function. To simplify the notation, put $z = x - b$, and expand it:

$$\frac{q}{2}\left[\frac{z}{q}\right]\left(\left[\frac{z}{q}\right]+1\right) + b\left[\frac{z}{q}\right] = \frac{q}{2}(z/q - \{z/q\})(z/q - \{z/q\} + 1) + b(z/q - \{z/q\})$$

$$= \frac{1}{2q}z^2 - z\{z/q\} + \frac{q}{2}\{z/q\}^2 + \frac{z}{2} - \frac{q}{2}\{z/q\} + b(z/q - \{z/q\}) \qquad (17)$$

$$= \frac{1}{2q}z^2 + o(\frac{1}{q}z^2).$$

Replacing $z = x - b$ back into the equation yields the stated asymptotic formula. Here, the value $x^2/2q + o(x^2/q) \geq 1$ is nonnegative for any admissible triples $0 < a < q \leq x$. Lastly, in (i), each root generates a distinct arithmetic progression, and there are $\rho(q) \geq 1$ distinct roots $b \in \{ x : x^3 \equiv a \bmod q \}$, so the maximal of the finite sum is the value $\rho(q)x^2/2q + o(\rho(q)x^2/q)$. ∎

The above estimate is sufficient for the intended application. A sharper estimate of the form $x^2/2q + O(x/q) + 1/4$ for (17) appears in [SP, p. 83].

**The Dirichlet Series**

The main interest here is in obtaining an effective estimate of the summatory function of the Dirichlet series

$$D(s,f) = \sum_{n \in D_f} \frac{\mu(n)\log n}{n^s} \tag{18}$$

of a complex variable $s \in \mathbb{C}$. This series is supported on the sub set of integers

$$\begin{aligned} D_f &= \{ d \in \mathbb{N} : x^3 + 2 \equiv 0 \bmod d \text{ for some } x \in \mathbb{Z} \} \\ &= \{ d \in \mathbb{N} : p \mid n \Rightarrow p = 3b+2, \text{ or } p \mid n \Rightarrow p = n = a^2 + 27b^2, \text{ or } p = 3 \text{ for } a, b \in \mathbb{Z} \}, \end{aligned} \tag{19}$$

prescribed by the polynomial $f(x) = x^3 + 2 \in \mathbb{Z}[x]$.

**Lemma 5.** The Dirichlet series $D(s,f) = \sum_{n \in D_f} (\mu(n)\log n) n^{-s}$ is analytic, and zerofree on the complex half plane $\{ s \in \mathbb{C} : \Re\mathrm{e}(s) \geq 1 \}$. In particular, the partial sum satisfies the relation

$$\sum_{n \leq x,\ n \in D_f} \frac{\mu(n)\log n}{n^s} = -\kappa_3 + O\left( x^{1-s} e^{-c(\log x)^{1/2}} \log x \right), \tag{20}$$

where $s \in \mathbb{C}$ is a complex number such that $\Re\mathrm{e}(s) \geq 1$, and $\kappa_3 > 0$, and $c > 0$ are constants.

*Proof*: The subsets of integers $D_f$ and $\mathbb{N} - D_f$ have positive densities $\delta(D_f) = \alpha > 0$, and $\delta(\mathbb{N} - D_f) = 1 - \alpha$ respectively, in the set of integers $\mathbb{N}$. Moreover, the derivative of the inverse zeta function $\zeta^{-1}(s) = \sum_{n \geq 1} \mu(n) n^{-s}$ is precisely the expression $\zeta'(s)/\zeta^2(s) = \sum_{n \geq 1} (\mu(n)\log n) n^{-s}$. From these data, it follows that

$$\begin{aligned} \sum_{n \in D_f} \frac{\mu(n)\log n}{n^s} &= \sum_{n \geq 1} \frac{\mu(n)\log n}{n^s} - \sum_{n \notin D_f} \frac{\mu(n)\log n}{n^s} \\ &= \frac{\zeta'(s)}{\zeta(s)^2} - (1-\alpha) \cdot \frac{\zeta'(s)}{\zeta(s)^2} \\ &= \alpha) \cdot \frac{\zeta'(s)}{\zeta(s)^2} \end{aligned} \tag{21}$$

is analytic and zerofree on the complex half plane $\{\, s \in \mathbb{C} : \mathfrak{Re}(s) \geq 1\,\}$. Hence, the summatory function satisfies

$$\sum_{n \leq x,\ n \in D_f} \frac{\mu(n) \log n}{n^s} = \alpha \cdot \frac{\zeta'(s)}{\zeta(s)^2} + O\!\left( x^{1-s} e^{-c(\log x)^{1/2}} \log x \right), \tag{22}$$

where $\kappa_3 = -\alpha \zeta'(s)\zeta(s)^{-2} > 0$, and $c > 0$ are constants. ∎

It appears that the Dirichlet series $D(s,f) = \sum_{n \in D_f} (\mu(n)\log n) n^{-s}$ can also be written in terms or approximated by a finite sum of the inverses $\zeta_Q^{-1}(s) = \sum_{n \geq 1} \mu(n) r_Q(n) n^{-s}$ of the Epstein zeta functions

$$\zeta_Q(s) = \sum_{(n,m) \neq (0,0)} \frac{1}{(am^2 + bmn + cn^2)^s} = \sum_{n \geq 1} \frac{r_Q(n)}{n^s}, \tag{23}$$

where $Q(x, y) = ax^2 + bxy + cy^2$, and $r_Q(n) = \{\, (x, y) : n = ax^2 + bxy + cy^2,\ x, y \in \mathbb{Z}\,\}$ with $\gcd(a, b, c) = 1$. Thus, the analysis of the Dirichlet series $D(s,f) = \sum_{n \in D_f} (\mu(n)\log n) n^{-s}$ can also be done using these standard analytic number theory tools.

**Estimate for Finite Sum over Prime Powers**
Finite sums of vonMangoldt function over the prime powers occur frequently in the analysis of prime numbers. An instance of these estimates is computed here.

***Lemma 6.*** Let $x \geq 1$ be a large real number. If $\Lambda(n^3 + k) \neq \Lambda(p)$ with $p$ prime for all $n^3 + k \leq x$, then

$$\sum_{n^3 + k = p^v \leq x,\ v \geq 2} n \Lambda(n^3 + k) \leq c x^{1/2} \log^2 x, \tag{24}$$

for some constant $c > 0$.

*Proof*: For a fixed $k$, there are at most $x^{1/3}$ integers of the form $n^3 + k \leq x$ in the interval $[1, x]$. Of these, a total of at most $(x^{1/3})^{1/2} \log x = x^{1/6} \log x$ are square integers $n^3 + k = m^2$, or powers $n^3 + k = m^v$, $v > 1$, and $m, n \in \mathbb{N}$. In light of these observations, the estimate is as follows:

$$\sum_{n^3 + k = p^v \leq x,\ v \geq 2} n \Lambda(n^3 + k) \leq x^{1/3} \sum_{p^v \leq x,\ v \geq 2} \log p^v \leq x^{1/3} \left( \frac{\log x}{\log 2} \sum_{n \leq x^{1/6}} \log n \right) \log x \leq c x^{1/2} \log^2 x, \tag{25}$$

where $v \leq \log(x)/\log(p) \leq \log(x)/\log(2)$, and $c > 0$ is a constant. ∎